\documentclass{amsart}

\usepackage{amsmath, amssymb, amsfonts, amsthm, amscd}
\usepackage{manfnt}
\usepackage{mathtools}
\usepackage[all]{xy}
   \SelectTips{cm}{10}
\usepackage{txfonts}
\newtheorem{thm}{Theorem}[section]
\newtheorem{lemma}[thm]{Lemma}
\newtheorem{prop}[thm]{Proposition}
\newtheorem{cor}[thm]{Corollary}

\newtheorem{prop-conj}[thm]{Proposition-Conjecture}

\theoremstyle{definition}
\newtheorem{defn}[thm]{Definition}

\theoremstyle{remark}
\newtheorem{rmk}[thm]{Remark}

\theoremstyle{remark}
\newtheorem{notation}[thm]{Notation}

\theoremstyle{remark}

\theoremstyle{remark}
\newtheorem{eg}[thm]{Example}

\newcommand{\Q}{\mathbb{Q}}
\newcommand{\Qb}{\overline{\mathbb{Q}}}
\newcommand{\Z}{\mathbb{Z}}
\newcommand{\CC}{\mathbb{C}}
\newcommand{\RR}{\mathbb{R}}

\DeclareMathOperator{\Hom}{Hom}
\DeclareMathOperator{\Ext}{Ext}
\DeclareMathOperator{\End}{End}
\DeclareMathOperator{\Aut}{Aut}

\DeclareMathOperator{\Ad}{Ad}

\DeclareMathOperator{\im}{im}

\DeclareMathOperator{\Res}{Res}

\newcommand{\gal}[1]{\Gamma_{#1}} 
\newcommand{\Gal}{\mathrm{Gal}} 

\newcommand{\into}{\hookrightarrow}
\newcommand{\onto}{\twoheadrightarrow}

\newcommand{\mc}{\mathcal}
\newcommand{\mf}{\mathfrak}
\newcommand{\mr}{\mathrm}

\newcommand{\mbb}{\mathbb}

\newcommand{\tT}{\widetilde{T}}
\newcommand{\tZ}{\widetilde{Z}}

\newcommand{\tM}{\widetilde{M}}

\begin{document}
\title{Mumford-Tate groups of polarizable Hodge structures}

\author{Stefan Patrikis}
\email{patrikis@math.mit.edu}
\address{Department of Mathematics\\MIT\\Cambridge, MA 02139}


\thanks{I am grateful to Phillip Griffiths for his interest, and for his feedback on the first draft of this note. This work was carried out with the support of NSF grant DMS-1062759, while I was a member of the Institute for Advanced Study in 2012-2013.}

\begin{abstract}
We classify the possible Mumford-Tate groups of polarizable rational Hodge structures.
\end{abstract}

\maketitle

\section{Introduction}
This note offers an answer to a question several authors (\cite{moonen:MTgroups} and \cite{griffiths-green-kerr:MTgroups}, for instance) have raised, that of describing which (connected, reductive) algebraic groups over $\Q$ can arise as Mumford-Tate groups of polarizable Hodge structures. The study \cite{griffiths-green-kerr:MTgroups} made much progress on this question, answering it for simple $\Q$-groups of adjoint type and all absolutely simple $\Q$-groups. We will simplify their arguments in a way that allows for a uniform treatment of the general case. 

The answer reached is very easy to work with, as we try to demonstrate through the examples of \S \ref{explaining}. A notable application-- of the techniques used in the proof, rather than of the classification itself-- is a (polarized) Hodge-theoretic analogue of a conjectural property of motivic Galois groups suggested by Serre (\cite[8.1]{serre:motivicgalois}); after posting the first version of this paper, I learned that Milne has proven essentially this same result in \cite[Corollary 8.6]{milne:shimvarandmoduli}.\footnote{Milne's result is stated only for Hodge structures satisfying Deligne's condition `SV1' that the Hodge co-character act with weights $-1, 0, 1$ on the Lie algebra of the Mumford-Tate group; but this is inessential to the argument.} The original motivic variant is closely related to a vast (conjectural) generalization of a classical construction of Kuga-Satake that associates to any complex $K3$ surface a complex abelian variety: see the notes \cite{stp:variationsarxiv}, which study arithmetic refinements and analogues of this problem.

The techniques of this paper are classical, building on ideas already present in \cite{griffiths-green-kerr:MTgroups}, \cite{serre:motivicgalois}, and \cite{deligne:weilK3}. One reason to be interested in the main theorem (\ref{main}) is that it provides the best-known (and, if one is an optimist, perhaps sharp) upper-bound on the possible collection of motivic Galois groups of pure motives over $\CC$.

\section{Setup and statement of the theorem} 
A useful reference for the basic notions described here is \cite{moonen:MTgroups}. Recall that a real Hodge structure is a representation $h \colon \mbb{S} \to \mr{GL}(V_{\RR})$ of the Deligne torus $\mbb{S}= \Res_{\CC/\RR}(\mbb{G}_m)$. A rational Hodge structure (abbreviated ``$\Q$-HS'') is a $\Q$-vector space $V$ along with a real Hodge structure on $V_\RR= V \otimes_{\Q} \RR$ such that the (weight) homomorphism $w_h= h|_{\mbb{G}_{m, \RR}}$ is defined over $\Q$. $V$ is pure of weight $n$ if $w_h(r)= r^{n} \cdot \mr{id}_V$.\footnote{This is the sign convention of \cite{griffiths-green-kerr:MTgroups}. The opposite convention is often used, especially in the theory of Shimura varieties.} Much more fundamental in algebraic geometry is the class of polarizable Hodge structures: a $\Q$-Hodge structure $V$ of weight $n$ is polarizable (or sometimes `$h(\imath)$-polarizable') if there exists a morphism of $\Q$-HS $Q \colon V \otimes V \to \Q(-n)$ with the critical positivity property that, on $V_{\RR} \otimes V_{\RR}$, the pairing
\[
(v, w) \mapsto (2 \pi \imath)^n Q(v, h(\imath)w)
\]
is symmetric and positive-definite. Define the category $\Q \mr{HS}^{pol}$ of (not necessarily pure) polarizable $\Q$-HS to have objects given by direct sums of pure polarizable $\Q$-HS's, and to have morphisms simply given by morphisms of $\Q$-HS's. With this definition, a basic observation is that $\Q\mr{HS}^{pol}$ is a semi-simple ($\Q$-linear) Tannakian category. We let $\mc{H}$ denote its Tannakian group, so that $\mc{H}$ is a connected pro-reductive group over $\Q$. There is a canonical (weight) homomorphism $w \colon \mbb{G}_{m, \Q} \to \mc{H}$; it maps to the center of $\mc{H}$.

A Mumford-Tate group-- for the purposes of this note, this will always mean Mumford-Tate group of a polarizable $\Q$-HS-- is simply any algebraic quotient of $\mc{H}$. More concretely, the Mumford-Tate group of a $\Q$-HS $V$ is the $\Q$-Zariski closure of the image of $h \colon \mbb{S} \to \mr{GL}(V_{\RR})$ (see, eg, Key Property $1.5$ of \cite{moonen:MTgroups}). We want to describe which (connected, reductive) $\Q$-groups arise in this way. Having articulated the main problem, to address it we must introduce certain other basic objects.

It is natural to study $\mc{H}$ in three steps:
\begin{enumerate}
\item Understand its center, or the center's neutral component $\mc{C}$, along with the (pro-)isogeny $\mc{C} \to \mc{S}$ onto the maximal abelian quotient $\mc{S}$ (the `Serre group') of $\mc{H}$.
\item Understand its derived group $\mc{D}$.
\item Understand how $\mc{C}$ and $\mc{D}$ are glued together under $\mc{C} \cdot \mc{D} = \mc{H}$.
\end{enumerate}
It is this third step that we have only partially achieved. Understanding the structure of $\mc{S}$ translates to the classical theory of Hodge structures with complex multiplication. The second step, for our purposes, amounts to the analogue of Serre's question (\cite[8.1]{serre:motivicgalois}) mentioned in the introduction, along with the description, due to Griffiths, Green, and Kerr (\cite{griffiths-green-kerr:MTgroups}) in the $\Q$-simple case, of which adjoint groups are realizable as Mumford-Tate groups. We will elaborate later on this background, but for now let us state the main result. To save breath, from now on we will simply say a group `is a M-T group' if it is the Mumford-Tate group of an object of $\Q\mr{HS}^{pol}$.

Some simple consequences of polarizability allow us to restrict the class of groups considered: 
\begin{notation}\label{hypotheses1}
Let $\tM$ be a connected reductive group over $\Q$. We enshrine the following notation:
\begin{itemize}
\item $M$ denotes the derived group of $\tM$, and $M^{ad}$ its adjoint group;
\item $\mf{m}$ is the Lie algebra of $M$;
\item If $\tM$ is the Mumford-Tate group of a $\Q$-HS $(V, h)$, then (by definition) $w= w_h$ is defined over $\Q$, and this easily implies it is central in $\tM$; if $V$ is polarizable, then (Remark $1.18$ of \cite{moonen:MTgroups}) $Z_{\tM}/w(\mbb{G}_{m})$ is compact over $\RR$, and we define a $\Q$-group $M'$ to be the connected normal complement to $w(\mbb{G}_m)$ in $\tM$: it descends $M'_{\RR}$, the subgroup generated by $M_{\RR}$ and the maximal compact sub-torus of $Z_{\tM}(\RR)$. 
\end{itemize}
\end{notation}
There is another obviously necessary condition to be a M-T group; this is the trivial reductive variant of IV.A.2 of \cite{griffiths-green-kerr:MTgroups}:
\begin{lemma}\label{necessary}
If $\tM$ is the M-T group of $(V, \tilde{h})$, then $M'_{\RR}$ contains a compact maximal torus $T'_{\RR}$ such that $\tilde{h}$ factors through the maximal torus $\tT_{\RR}= w(\mbb{G}_{m, \RR})\cdot T'_{\RR}$ of $\tM_{\RR}$.
\proof
We may assume $V$ is pure, so that there is a polarizing form $Q$ and a faithful representation $\rho \colon \tM \into G(V, Q)$, where $G(V, Q)$ denotes the similitude group for the pairing $Q$ (similarly, we will denote by $U(V, Q)$ the automorphism group of the pairing). Now consider a maximal torus $T'_{\RR}$ of $M'_{\RR}$ that contains $\tilde{h}(\mbb{S}^1)$ (itself clearly contained in $M'(\RR)$. Then $T'_{\RR}$ is contained in the centralizer of $h(\mbb{S}^1)$, which is a closed subgroup of the isotropy group-- necessarily compact-- of $h|_{\mbb{S}^1}$ in $U(V, Q)_{\RR}$. Thus $T'(\RR)$ is compact.
\endproof
\end{lemma}
\begin{notation}\label{hypotheses2}
From now on, $\tM$ will denote a connected reductive group over $\Q$ satisfying the hypotheses:
\begin{itemize}
\item There exists a central homomorphism $w \colon \mbb{G}_{m, \Q} \to \tM$ such that $Z_{\tM}/w(\mbb{G}_m)$ is compact over $\RR$.\footnote{Note that $w$ is not `part of the data'-- it could be replaced by some non-zero power, for instance. What is canonical is the subgroup $w(\mbb{G}_{m, \Q})$.}
\item $M^{ad}_{\RR}$ contains a compact maximal torus; we will typically denote one by $T^{ad}_{\RR}$, with $T^{ad}$ denoting the base-change to $\CC$.
\end{itemize}
\end{notation}
Griffiths, Green, and Kerr have described which $h \colon \mbb{S}^1 \to T^{ad}_{\RR} \subset M^{ad}_{\RR}$ can lead to polarized Hodge structures on the semisimple Lie algebra $\mf{m}$. We will refer to (a slight variant-- see Remark \ref{cancochar}-- of) the set of such $h$ as the `polarizable congruence classes' of co-characters in $X_{\bullet}(T^{ad})$,\footnote{Note that there is a canonical bijection between homomorphisms from $\mbb{S}^1$ to the compact torus $T^{ad}_{\RR}$ and co-characters $\mbb{G}_{m, \CC} \to T^{ad}$ of the complexified torus $T^{ad}:= T^{ad}_{\RR} \otimes \CC$.} which we now describe. Let $\theta$ be a Cartan involution of $\mf{m}_{\RR}$ (recall that all such $\theta$ are $M^{ad}(\RR)$-conjugate), yielding a Cartan decomposition $\mf{m}_{\RR}= \mf{k}_{\RR} \oplus \mf{p}_{\RR}$ into $\theta=1$ and $\theta=-1$ eigenspaces, and let $\mf{t}_{\RR}= \mr{Lie}(T^{ad}_{\RR}) \subset \mf{k}_{\RR}$ be a compact maximal torus inside $\mf{k}_{\RR}$. This yields a partition of the roots of $\mf{t}= \mf{t}_{\RR} \otimes \CC$ in $\mf{m}_{\CC}$ into compact (contained in $\mf{k}$) and non-compact (contained in $\mf{p}$) roots.
\begin{defn}[compare Proposition IV.B.3 of \cite{griffiths-green-kerr:MTgroups}]
With this fixed Cartan involution $\theta$, and choice of a $\theta$-stable compact maximal torus $T^{ad}_{\RR}$, the polarizable congruence classes in $X_{\bullet}(T^{ad})$ consist of those co-characters $\mu$ (`$\theta$-polarizable') such that
\begin{align*}
& \text{$\langle \mu, \alpha \rangle$ is even for all compact roots $\alpha$;}\\
& \text{$\langle \mu, \beta \rangle$ is odd for all non-compact roots $\beta$.}
\end{align*}
\end{defn}
\begin{rmk}\label{cancochar}
\begin{enumerate}
\item The authors of \cite{griffiths-green-kerr:MTgroups} consider instead the co-characters $\mf{l} \colon \mbb{S}^1 \to T^{ad}_{\RR}$ that satisfy the related congruences 
\begin{align*}
&\text{$\langle \mf{l}, \alpha \rangle \equiv 0 \pmod 4$ for all compact roots $\alpha$;} \\
&\text{$\langle \mf{l}, \beta \rangle \equiv 2 \pmod 4$ for all non-compact roots $\beta$.}
\end{align*}
When such an $\mf{l}$ is extended (by declaring it trivial on $\mbb{G}_{m, \RR}$) to a homomorphism $h \colon \mbb{S} \to M^{ad}_{\RR}$, what we denote by $\mu$ is simply the associated Hodge co-character, and is easily checked to equal $\frac{1}{2} \mf{l}$ in $X_{\bullet}(T^{ad})$: namely, $\Ad \mf{l}(z)$ acts on $\mf{m}^{-j, j}$ by $z^{-j} \bar{z}^{j}= z^{-2j}$, while $\Ad \mu(z)$ by definition acts as $z^{-j}$.
\item This property does not depend on the choice of the $\theta$-fixed compact torus $T^{ad}_{\RR}$ containing the image of $\mf{l}$: that is, given a $\mu$ whose corresponding $\mf{l} \colon \mbb{S}^1 \to T^{ad}_{\RR}$ also factors through a second $\theta$-fixed maximal torus $S_{\RR}^{ad}$, $(M^{ad}_{\RR})^{\theta=1}$-conjugacy of the two tori implies $\theta$-polarizability of $\mf{l}$ is independent of whether we decompose $\mf{m}_{\CC}$ with respect to $T^{ad}$ or $S^{ad}$.
Note too that if $\mu$ is polarizable with respect to $\theta$, then for all $x \in M^{ad}(\RR)$, $\Ad(x) \mu$ is polarizable with respect to $\Ad(x) \theta$.
\item The relations $[\mf{k}_{\RR}, \mf{k}_{\RR}]\subset \mf{k}_{\RR}$, $[\mf{k}_{\RR}, \mf{p}_{\RR}]\subset \mf{p}_{\RR}$, and $[\mf{p}_{\RR}, \mf{p}_{\RR}]\subset \mf{k}_{\RR}$ imply that it is equivalent to check the above congruences on just the simple (for some choice of dominant Weyl chamber) compact and non-compact roots. The full preimage in $X_{\bullet}(T^{ad})$ of some subset of $X_{\bullet}(T^{ad})/2X_{\bullet}(T^{ad})$ gives all such $\mu$.
\end{enumerate}
\end{rmk}
The theorem we will now state contains (and generalizes to the reductive case) the main results on Hodge representations in \S IV.B and \S IV.E of \cite{griffiths-green-kerr:MTgroups}. The argument is also much simpler: essentially, we avoid that paper's (`main theorem on Hodge representations') Theorem IV.B.6 (and its intricate sub-result IV.E.4) by proving, directly and more generally, its Corollary IV.B.8; this latter point follows from the observations in \S $2$ of \cite{deligne:weilK3}.
\begin{thm}\label{main}
Let $\tM$ satisfy the hypotheses of Notation \ref{hypotheses2}, which we have seen are necessary for $\tM$ to be a Mumford-Tate group. Then $\tM$ is in fact the Mumford-Tate group of some polarizable $\Q$-Hodge structure if, and only if, the following two conditions are satisfied:
\begin{itemize}
\item For some (any) Cartan involution $\theta$ of $M^{ad}_{\RR}$ and some (any) $\theta$-stable maximal compact torus $T^{ad}_{\RR}$ of $M^{ad}_{\RR}$, some member of the polarizable congruence classes in $X_{\bullet}(T^{ad})$ lifts to $\tM_{\CC}$.
\item The neutral component $Z^0_{\tM}$ of the center of $\tM$ is a quotient of $\mc{C}$, or, equivalently, of $\mc{S}$.
\end{itemize}
\end{thm}
The condition on the center just says that $Z^0_{\tM}$ is the M-T group of a CM Hodge-structure. Also note that it is possible for some, but not all, members of the polarizable congruence classes to lift to $\tM_{\CC}$.
\section{Explaining the theorem}\label{explaining}
In this section we will give a few examples and corollaries of Theorem \ref{main}. First we re-package the statement that some member of the polarizable congruence classes in $X_{\bullet}(T^{ad})$ lifts to $\tM_{\CC}$. Let $\tT$ denote the pre-image of $T^{ad}$ in $\tM_{\CC}$, so that we have an exact sequence of groups of multiplicative type:
\[
1 \to Z_{\tM, \CC} \to \tT \to T^{ad} \to 1.
\]
Taking character groups and applying $\Hom(\cdot, \Z)$, we get an exact sequence of abelian groups:
\[
\xymatrix{
0 \ar[r] & X_{\bullet}(Z_{\tM, \CC}) \ar[r] & X_{\bullet}(\tT) \ar[r] & X_{\bullet}(T^{ad}) \ar[r] \ar@{-->}[d] & \Ext^1(X^{\bullet}(Z_{\tM, \CC}), \Z) \ar[d]^{\sim} \\
 & & & \left(X^{\bullet}(Z_{\tM, \CC})_{tor}\right)^D & \ar[l]_{\sim} \Ext^1(X^{\bullet}(Z_{\tM, \CC})_{tor}, \Z), \\
 }
\]
where $A^D$ denotes $\Hom(A, \Q/\Z)$ for an abelian group $A$. A class $\mu \in X_{\bullet}(T^{ad})$ lifts to $\tM_{\CC}$ if and only if its image under the boundary map, or, equivalently, its image in $\left(X^{\bullet}(Z_{\tM, \CC})_{tor}\right)^D$, is zero. We conclude:
\begin{cor}\label{notwotorsion}
The condition that some member of the polarizable congruence classes of $X_{\bullet}(T^{ad})$ lift to $\tM_{\CC}$ amounts to asking that for one of its elements $\mu$, the image of $\mu$ in
\[
\left(X^{\bullet}(Z_{\tM, \CC})_{tor}\right)^D/2\left(X^{\bullet}(Z_{\tM, \CC})_{tor}\right)^D \cong \left( X^\bullet(Z_{\tM, \CC})[2] \right)^D
\]
be zero. In particular, if $X^\bullet(Z_{\tM, \CC})$ has no two-torsion-- in particular, if $Z_{\tM}$ is a torus!-- then there is no lifting obstruction, so that $\tM$ is a M-T group provided $Z^0_{\tM}$ is a quotient of $\mc{S}$.
\end{cor}
In the next result, we use arguments similar in spirit to those of \S \ref{proof} to obtain the desired application to (the polarized Hodge-theoretic version of) Serre's question. This result is originally due to Milne, with a slightly different argument, in \cite[Propositions 8.1, 8.5, Corollary 8.6]{milne:shimvarandmoduli}:
\begin{cor}
As before, let $\mc{H}$ denote the Tannakian group of $\Q\mr{HS}^{pol}$, and let $\mc{D}$ denote its derived group. Then $\mc{D}$ is simply-connected.
\end{cor}
\proof
We must show that for any surjection $\rho_{\mc{D}} \colon \mc{D} \onto H$, with $H$ a semi-simple group over $\Q$, there is a lift $\mc{D} \to H^{sc}$ to the simply-connected cover of $H$. Let us write $H^{sc} \xrightarrow{\alpha} H \xrightarrow{\beta} H^{ad}$ for the covering maps. The given $\rho_{\mc{D}}$ induces a surjection
\[
\rho \colon \mc{H} \onto \mc{D}/Z_{\mc{D}} \onto H^{ad}.
\]
We will embed $Z(H^{sc})$ into a well-chosen $\Q$-torus $\tZ$ that will allow us to lift $\mc{H} \onto H^{ad}$ to a homomorphism 
\[
\mc{H} \to \widetilde{H}= H^{sc} \underset{Z(H^{sc})}{\times} \widetilde{Z}.
\]
We first claim that the $\gal{\Q}= \Gal(\Qb/\Q)$-module $X^\bullet(Z(H^{sc}))$ splits over the maximal CM extension $\Q^{cm}$ of $\Q$, and that (any) complex conjugation $c \in \gal{\Q}$ acts as $-1$ on $X^\bullet(Z(H^{sc}))$. The first claim holds provided some maximal torus $T^{sc}$ of $H^{sc}$ splits over $\Q^{cm}$, hence provided some maximal torus $T^{ad}$ of $H^{ad}$ splits over $\Q^{cm}$. But $H^{ad}$ is a Mumford-Tate group, and the existence of polarizations formally implies that $H^{ad}$ (or indeed $\mc{H}$) is split over $\Q^{cm}$ (compare \cite[Lemma 16.3.1]{stp:variationsarxiv}). For the second claim, let $T^{ad}_{\RR}$ be a compact (see Lemma \ref{necessary}) maximal torus in $H^{ad}_{\RR}$ through which the canonical composition $h \colon \mathbb{S} \to \mc{H}_{\RR} \to H^{ad}_{\RR}$ factors,\footnote{The notation should \textit{not} be misread as indicating that $T^{ad}_{\RR}$ arises by base-change from a torus over $\Q$: this is the case if and only if $\rho$ is a CM Hodge structure.} and let $T^{sc}_{\RR}$ be the pre-image in $H^{sc}_{\RR}$. Since $Z(H^{sc})_{\RR}= Z(H^{sc}_{\RR})$ is a subgroup of the compact torus $T^{sc}_{\RR}$, $c$ acts as $-1$, as claimed, on the character group. Thus, not only can we embed $Z(H^{sc})$ into a torus of the form $\Res_{K/\Q}(\mathbf{G}_m)^k$ for some CM field $K$ and some integer $k$, but also the quotient map
\[
X^{\bullet}\left(\Res_{K/\Q}(\mathbf{G}_m)^k \right) \onto X^\bullet \left( Z(H^{sc}) \right)
\]
factors through the maximal quotient on which $c$ acts by $-1$. Otherwise put, there is an embedding $Z(H^{sc}) \into U_K^k$, where $U_K$ is the $\Q$-torus with character group
\[
X^\bullet(U_K)= \left( \bigoplus_{\tau \colon K \into \CC} \Z \tau \right)/ \sum_{\tau} \Z(\tau + c \tau).
\]
We then let $\tZ= U_K^k$, and we enlarge $H^{sc}$ to the group $\widetilde{H}= H^{sc} \underset{Z(H^{sc})}{\times} \widetilde{Z}$, whose adjoint group is canonically isomorphic to $H^{ad}$, but whose center is now a torus. Now, the restriction of $h$ to $\mathbf{G}_{m, \RR}$ is trivial since $H^{ad}$ is adjoint, hence $h$ factors through $\mathbb{S}/\mathbf{G}_{m, \RR} \cong \mathbb{S}^1$. This $\mathbb{S}^1$-representation lifts to $\widetilde{H}_{\RR}$: the obstruction to lifting lies in
$\Ext^1_{\Z[\gal{\RR}]}(X^\bullet(\tZ), X^\bullet(\mathbb{S}^1)$,
but since $\tZ_{\RR}$ is just a number of copies of $\mathbb{S}^1$, it is easy to see that this group of extensions is trivial. In particular, we have a weight-zero (and thus having weight homomorphism defined over $\Q$) lift $\tilde{h} \colon \mathbb{S} \to \widetilde{H}_{\RR}$ of $h$, and we deduce that $\tilde{h}$ is polarizable, since the center $\tZ(\RR)$ is compact (see Lemma \ref{deligne}). Hence $\tilde{h}$ arises from a $\tilde{\rho} \colon \mc{H} \to \widetilde{H}$, necessarily lifting $\rho$. Restricting to $\mc{D}$, we obtain a surjection $\tilde{\rho} \colon \mc{D} \onto H^{sc}$ lifting $\beta \circ \rho_{\mc{D}}$. This means that $\alpha \circ \tilde{\rho}$ and $\rho_{\mc{D}}$ differ by a homomorphism $\mc{D} \to Z(H)$; but $\mc{D}$ is connected, so $\tilde{\rho}$ in fact lifts $\rho_{\mc{D}}$, completing the proof.
\endproof

We next want to clarify the relationship between our results and those of \cite{griffiths-green-kerr:MTgroups}. The important notational difference is that what they call a `Mumford-Tate group' (throughout Chapter IV of their monograph) is what is often called the Hodge group. Thus $\mr{SL}_2$ is a Mumford-Tate group in their sense, but not in ours; in their language, this is reflected in the statement that `the only faithful Hodge representations of $\mr{SL}_2$ are of odd weight.' In general, their problem of deciding whether one has faithful Hodge representations of even or odd weight is just the question of whether elements of the polarizable congruence classes lift to the group $M'$ (even weight) or require a copy of $\mbb{G}_{m, \Q}$ inside $\tM$ in order to lift. In our language, to say whether $M$ (now assumed semi-simple) is a Hodge group we just replace the lifting criterion for elements $\mu$ of the polarizable congruence classes with the corresponding lifting criterion for $\mf{l}= 2\mu$; in particular, this recovers the main Theorem IV.B.6 of \cite{griffiths-green-kerr:MTgroups}.\footnote{Doubling our $\pmod 2$ condition gives their $\pmod 4$ condition; note that this makes it easier to lift, because, for instance, $X^{\bullet}(Z_{\tM, \CC})$ must now have \textit{four}-torsion-- not merely two-torsion-- for there to be a lifting obstruction.} We will not dwell further on this bookkeeping.

Now we give an example to show how to determine liftability of co-characters $\mu \in X_{\bullet}(T^{ad})$. This is in all cases an easy root/weight calculation, which is carried out for absolutely simple $\Q$-groups in \cite[IV.E]{griffiths-green-kerr:MTgroups}; note that the lifting criterion in Theorem \ref{main} only depends on the underlying real group $\tM_{\RR}$, so we can ignore the $\Q$-structure and perform the same sorts of root/weight calculations in general. 
\begin{eg}
We work out the relatively interesting example of even orthogonal groups, i.e. the various real forms of $\mf{m}_{\CC}= \mf{so}_{2n}(\CC)$. In (usual) coordinates, we can write $X_{\bullet}(T^{ad})= \oplus_{i=1}^{n} \Z \lambda_i + \Z(\frac{\sum \lambda_i}{2})$. $X^{\bullet}(T^{ad})$ is then the subspace of $\oplus_{i=1}^{n} \Z \chi_i$ of elements $\sum c_i \chi_i$ satisfying $\sum c_i \equiv 0 \pmod 2$. Here $\langle \chi_i, \lambda_j \rangle= \delta_{i, j}$, and for simple roots we take 
\[
\alpha_1= \chi_1 - \chi_2, \ldots, \alpha_{n-2}= \chi_{n-2}- \chi_{n-1}, \alpha_{n-1}= \chi_{n-1}- \chi_{n}, \alpha_{n}= \chi_{n-1}+ \chi_{n}.
\]
The classification of (absolutely) simple real Lie algebras (see \cite[\S VI.10]{knapp:beyond}) implies that any real form of $\mf{so}_{2n}(\CC)$ associated to a Vogan diagram (\cite[VI.8]{knapp:beyond}) with trivial $\theta$-action is given by a Vogan diagram in which the simple roots appear ordered as we have just ordered them, and with exactly one painted (i.e., non-compact) simple root: for $1 \leq p \leq n-2= p+q-2$, painting $\alpha_p$ gives the form $\mf{so}(2p, 2q)$; painting either $\alpha_{p+q-1}$ or $\alpha_{p+q}$, we obtain the real form $\mf{so}^*(2p+2q)$.

Assume $q \geq 2$, so we are in the $\mf{so}(2p,2q)$ case. Writing $\mu= \sum a_i \lambda_i$ with all $a_i \in \frac{1}{2} \Z$ and congruent mod $\Z$, the condition to lie in the polarizable congruence classes is that all $a_i$ be integers with 
\begin{align*}
&a_1 \equiv a_2 \equiv \ldots \equiv a_p \pmod 2 \\ 
&a_p \not \equiv a_{p+1} \pmod 2 \\
&a_{p+1} \equiv a_{p+2} \equiv \ldots \equiv a_{p+q} \pmod 2.
\end{align*}
In conclusion, all members of the polarizable congruence classes lift to $\mr{SO}(2p, 2q)$,\footnote{More precisely, the co-characters lift to $\mr{SO}_{2n}(\CC)$; the Hodge structures $h \colon \mbb{S} \to \mr{SO}(2p,2q)/\pm 1$ lift to $\mr{SO}(2p,2q)$.} and some member lifts to $\mr{Spin}(2p, 2q)$ precisely when either $p$ or $q$ is even. For the compact form $\mf{so}(2n)$, some member of the polarizable congruence classes will lift to $\mr{Spin}(2n)$; all members do precisely when $n$ is even. We stress, however, that in all cases the Hodge structures arising from any member of the polarizable congruence classes will lift to the appropriate form of $\mr{GSpin}_{2n}(\CC)$, since they lift to $\mr{SO}(2p,2q)$ and the kernel of $\mr{GSpin}(2p,2q) \to \mr{SO}(2p,2q)$ is $\mbb{G}_{m, \RR}$. 

When the painted root is either $\alpha_{n-1}$ or $\alpha_n$-- we may assume it is $\alpha_n$-- then we find $a_1 \equiv \ldots \equiv a_n \pmod 2$ and $a_{n-1}+a_n \equiv 1 \pmod 2$. This forces the $a_i$ all to lie in $\frac{1}{2}+\Z$, so there is no lift to $\mr{SO}^*(2n)$. When $n$ is odd there is no other lifting possibility, while for $n$ even there are two other double covers of $\mr{SO}^*(2n)/\pm1$ to which one might lift.

In the case where the Vogan diagram has a non-trivial $\theta$-orbit (consisting of $\alpha_{n-1}$ and $\alpha_n$), it is likewise known that any form is represented by a diagram with $\alpha_p$ ($1 \leq p \leq n-2=p+q-2$) painted, corresponding to $\mf{so}(2p+1, 2q-1)$, or with no simple root painted, corresponding to $\mf{so}(1, 2n-1)$. The same analysis as above would characterize when the polarizable conjugacy classes lift to $\mr{Spin}_{2n}(\CC)$, but $\mr{SO}(2p+1, 2q-1)$ does \textit{not} have a compact maximal torus, so no M-T group can have real Lie algebra $\mf{so}(2p+1, 2q-1)$.\footnote{This important point is missed in \cite[IV.E]{griffiths-green-kerr:MTgroups}.}

To give a geometric example, let $X/\CC$ be a \textit{projective} $K3$ surface. $H^2(X, \Q)(1)$ is a polarizable Hodge structure with M-T group $M$ embedded in a $\Q$-form of $\mr{SO}_{22}(\CC)$ whose real points are $\mr{SO}(3, 19)$; note that  $M_{\RR}$ cannot equal the full $\mr{SO}(3, 19)$, by the previous paragraph-- this reflects the existence of an ample line bundle (when $X$ is projective), which forces $M_{\RR}$ to embed into (and generically be equal to) $\mr{SO}(2,19)$. Since $h^{2,0}(X)=1$, the Hodge co-character is in this case $\mu= \lambda_1$ (which is indeed in the polarizable congruence class); it does not lift to $\mr{Spin}_{22}(\CC)$ (or $\mr{Spin}_{21}(\CC)$), but it does to $\mr{GSpin}_{22}(\CC)$ (or $\mr{GSpin}_{21}(\CC)$), and normalizing the lift to weight one, we obtain (via the action of $\mr{GSpin}$ on the Clifford algebra) the Kuga-Satake abelian variety associated to $X$. 

In contrast to the projective case, we have:
\begin{lemma}
Let $X$ be a sufficiently generic non-projective $K3$ surface; then the M-T group $M$ of the (non-polarizable!) Hodge structure $H^2(X, \Q)(1)$ will satisfy $M_{\RR}= \mr{SO}(3, 19)$. Moreover, if $X$ satisfies $\mr{Pic}(X)=0$ and $\End_{M}(H^2(X))= \Q$,\footnote{Assuming $\mr{Pic}(X)=0$, the transcendental `lattice' (except we work rationally) $T(X)$, defined to be the smallest $\Q$-Hodge sub-structure $T \subset H^2(X, \Q)(1)$ such that $T_{\CC} \supset H^{2,0}(X)(1)$, is equal to all of $H^2(X)$, i.e. $H^2(X)$ is an irreducible $\Q$-Hodge structure. Since $h^{2,0}=1$, $\End_M(H^2(X))$ is isomorphic to a subfield of $\CC$.} then $M_{\RR}= \mr{SO}(3, 19)$.
\end{lemma}
\proof
For the first claim, consider the holomorphically-varying family of Hodge structures over the entire $K3$ period domain.\footnote{I.e., letting $\Lambda$ denote the $K3$ lattice, the period domain is $D= \{x \in \mbb{P}(\Lambda_{\CC}): x \cdot x=0, x \cdot \bar{x}>0\}$.} The generic Mumford-Tate group (see the proof of Proposition \ref{simpleadjoint}) is a normal subgroup of the $\Q$-simple group $\mr{SO}(\Lambda_{\Q})$, hence they are equal. 

For the second claim, recall the theorem of Zarhin (\cite[2.2.1,2.3.1]{zarhin:K3hodge}) that describes the M-T group of an irreducible \textit{polarizable} Hodge structure $T$ of $K3$-type, showing it is determined by the field $E= \End_{MT_T}(T)$. In the non-polarizable case, we cannot deduce that $E$ is totally real or CM, so Zarhin's general result breaks down. But when $E=\Q$, one can check that the argument of \cite[\S 2.5]{zarhin:K3hodge} goes through unhindered. We omit the details.
\endproof
\end{eg}
We conclude with an example illustrating the condition on the center of a M-T group. Having a polarization implies (see \cite[1.23]{moonen:MTgroups}) that $Z_{M'}$ embeds in a product of groups of the form $U_F$, where $F$ is a CM or totally real field, and $U_F$ is the multiplicative $\Q$-group with $\Q$-points $\{x \in F^\times: x \bar{x}= 1\}$. The condition that $\mc{S}$ surject onto $Z^0_{\tM}$ is stronger, as seen from the following standard description of $X^\bullet(\mc{S})$ (\cite[7.3]{serre:motivicgalois}): letting $\Q^{cm}$ denote the maximal CM extension of $\Q$, with complex conjugation $c$, $X^{\bullet}(\mc{S})$ is 
\[
\{
\text{$f \colon \Gal(\Q^{cm}/\Q) \to \Z$: there exists an integer $w(f)$ satisfying $f(s)+f(sc)= w(f)$ for all $s$}
\}.
\]
The $\Gal(\Qb/\Q)$-action (factoring through $\Gal(\Q^{cm}/\Q)$) is $(t \cdot f)(s)= f(st)$. So, for example, letting $F/\Q$ be any quadratic imaginary extension, $U_F=U_F^0$ is a possible M-T group, but $U_F \times U_F$ is not. $X^\bullet(U_F)$ is isomorphic to a free $\Z$-module of rank $1$ on which $\Gal(\Qb/F)$ acts trivially, and on which $c$ acts by $-1$. Under any embedding $X^{\bullet}(U_F) \into X^\bullet(\mc{S})$, a generator must map to a function $f$ factoring through $\Gal(F/\Q)$ and satisfying $f(c)=-f(1)$; up to scaling, there is only one such $f$. One can generate many such examples, using the fact that $X^\bullet(\mc{S})$ is contained in the `regular representation' of $\Gal(\Q^{cm}/\Q)$ (so at least after tensoring with $\Q$, multiplicities of irreducible constituents are bounded by their ranks).
\section{Proof of the theorem}\label{proof}
We now give the proof of Theorem \ref{main}. One direction is simpler:
\begin{prop}\label{easydirection}
Suppose that $\tM$ is a M-T group. Recall that $\tM$ satisfies the conditions of Notation \ref{hypotheses2}. Moreover, there is a Cartan involution of $M^{ad}_{\RR}$ such that some member $\mu$ of the polarizable congruence classes in $X_{\bullet}(T^{ad})$ lifts to $\tM_{\CC}$, and $Z^0_{\tM}$ is a quotient of the Serre group $\mc{S}$.
\end{prop}
\proof
By assumption, there is a morphism $\tilde{h} \colon \mbb{S} \to \tM_{\RR}$, with $\Q$-Zariski dense image, and a faithful family $(V_i, Q_i)$ of representations $\rho_i \colon \tM \to G(V_i, Q_i)$ such that $Q_i$ polarizes the Hodge structure $\rho_i \circ \tilde{h}$ (in particular, $Q_i$ is $M'$-invariant). We will check that $\Ad \tilde{h}(\imath)$ is a Cartan involution of $M'(\RR)$. Assuming this for the time being, we consider the composition $h \colon \mbb{S} \to M^{ad}_{\RR}$ of $\tilde{h}$ with $\tM_{\RR} \onto M^{ad}_{\RR}$; $\Ad h(\imath)$ is then a Cartan involution of $M^{ad}(\RR)$. On the $\gamma$-root space in $\mf{m}_{\CC}$, $\Ad h(\imath)$ acts as $\imath^{\langle 2\mu_h, \gamma \rangle}$, where $\mu_h \in X_{\bullet}(T^{ad})$ is the Hodge co-character. Letting $\mf{m}_{\RR}= \mf{k}_{\RR} \oplus \mf{p}_{\RR}$ be the Cartan decomposition with respect to $\Ad h(\imath)$, it follows that $\mu_h$ lies in the corresponding polarizable congruence classes.

Now we check that $\Ad \tilde{h}(\imath)$ is indeed a Cartan involution of $M'_{\RR}$; we merely imitate \S  2 of \cite{deligne:weilK3} (and see Lemma \ref{deligne} below). By assumption, we have a morphism of Hodge-structures $Q_i \colon V_i \otimes V_i \to \Q(-n_i)$ (for some integers $n_i$). Let $(M')^{\tilde{\sigma}}$ denote the invariants in $M'(\CC)$ under $x \mapsto \Ad \tilde{h}(\bar{x})$ (where $\bar{x}$ is conjugation with respect to the real structure $M'_{\RR}$). Then the image of $(M')^{\tilde{\sigma}}$ in $\mr{GL}(V_{i, \CC})$ preserves the positive-definite hermitian form 
\[
(v, w) \mapsto (2 \pi \imath)^{n_i} Q_i(v, \tilde{h}(\imath) \bar{w}),
\]
hence is compact. Since $(V_i)$ is a faithful family, $(M')^{\tilde{\sigma}} \to \prod \mr{GL}(V_{i, \CC})$ is faithful, and thus $(M')^{\tilde{\sigma}}$ is compact, as claimed.

The claim about the center is easy: we have a surjection $\mc{H} \onto \tM$, hence a surjection $\mc{C} \onto Z^0_{\tM}$, hence (taking some power) also a surjection $\mc{S} \onto Z^0_{\tM}$.
\endproof
For the converse direction of Theorem \ref{main}, we proceed in three steps: simple adjoint groups, semi-simple adjoint groups, and finally the generally case. First we recall the basic lemma of Deligne (\cite[2.11]{deligne:weilK3}), slightly reformulated, that allows us to circumvent many of the more intricate arguments of \cite{griffiths-green-kerr:MTgroups}. 
\begin{lemma}[\cite{deligne:weilK3}]\label{deligne}
Let $\tM$ be as in Notation \ref{hypotheses2}, and suppose moreover that the rational (central) weight homomorphism $w \colon \mbb{G}_{m, \Q} \to \tM$ underlies the restriction $\tilde{h}|_{\mbb{G}_{m, \RR}}$ of a homomorphism $\tilde{h} \colon \mbb{S} \to \tM_{\RR}$. Then the following are equivalent:
\begin{enumerate}
\item $\Ad \tilde{h}(\imath)$ is a Cartan involution of $M'_{\RR}$;
\item Every homogeneous representation\footnote{We say $\rho \colon \tM \to \mr{GL}(V)$ is homogeneous of weight $n$ if $\rho \circ w(z)$ acts as $z^n$.} of $\tM$ is polarizable; 
\item $\tM$ has a faithful family of homogeneous, polarizable representations.
\end{enumerate}
\end{lemma}
\proof
$(3) \implies (1)$ was included in the proof of Proposition \ref{easydirection}. We now sketch the implications $(1)\implies(2)\implies(3)$. Let $V$ be a homogeneous representation of $\tM$. There are two steps: polarizability of $V_{\RR}$ comes from the fact that any representation of a compact Lie group (in this case the invariants $(M')^{\tilde{\sigma}}$ under $\tilde{\sigma}(x)= \Ad \tilde{h}(\imath)(\bar{x})$ inside $M'(\CC)$) has, by integration, an invariant inner product. This yields an $M'_{\RR}$-invariant $\tilde{h}(\imath)$-polarization on $V_{\RR}$, which we interpret as a $\tM_{\RR}$-invariant polarization $Q \colon V_{\RR} \otimes V_{\RR} \to \RR(-n)$, where by definition $w(x)$ acts on $\RR(-n)$ by $x^{2n}$.\footnote{That is, $Q(w(x)m\cdot v,w(x)m \cdot w)= x^{2n}Q(mv, mw)= x^{2n}Q(v, w)$.} To deduce the existence of a polarizing pairing on the $\Q$-vector space $V$, we note that (\cite[Lemme 2.10]{deligne:weilK3}), starting from the \textit{rational} vector space $P_{\Q}$ of $(-1)^{n}$-symmetric, $M'$-invariant, bilinear forms, the space of $\tilde{h}(\imath)$-polarizations in $P_{\RR}$ is open and non-empty, hence contains rational points.

Having shown $(2)$, we note that $\tM$ has a faithful family of homogeneous, hence polarizable, representations  because any irreducible (over $\Q$) representation $\tM \to \mr{GL}(V)$ must be homogeneous: $\End_{\Q[\tM]}(V)$ is a division algebra $D$ over $\Q$, and any homomorphism $\mbb{G}_{m, \Q} \to D^\times$ factors through a maximal torus, necessarily of the form $F^\times$ for some maximal commutative subfield $F \subset D$, and hence must be of the form $\Q^\times \xrightarrow{z \mapsto z^n} \Q^\times \subset F^\times \subset D^\times$.\endproof
The initial case in which $M$ is a $\Q$-simple adjoint group (so the hypotheses in Theorem \ref{main} are trivially satisfied) was established in \cite{griffiths-green-kerr:MTgroups}:
\begin{prop}[Theorem IV.E.1 of \cite{griffiths-green-kerr:MTgroups}]\label{simpleadjoint}
A $\Q$-simple adjoint group $M$ is a Mumford-Tate group if and only if $M_{\RR}$ contains a compact maximal torus.
\end{prop}
\begin{rmk}
The theorem is in fact misstated in \cite{griffiths-green-kerr:MTgroups}: the authors claim $M$ is a M-T group if and only if $M$ (over $\Q$) contains an anisotropic maximal torus, asserting that this is equivalent to $M_{\RR}$ having a compact maximal torus. This is false: let $M= \mr{SL}_3/\Q$ (or for the adjoint case $\mr{PGL}_3$), and choose a totally real cubic extension $F/\Q$. Then we get an embedding of the (two-dimensional, so maximal) norm-one torus $\Res_{F/\Q}^{1} \mbb{G}_m \into \mr{SL}_3/\Q$; thus $\mr{SL}_3/\Q$ contains an anisotropic maximal torus, but it is not a Mumford-Tate group, nor does $\mr{SL}_3(\RR)$ contain a compact maximal torus.
\end{rmk}
\proof
We review the proof, recasting slightly the argument in \cite[IV.A.9, IV.B.3, IV.E.1]{griffiths-green-kerr:MTgroups}. We need only treat the `if' direction, so assume that $M_{\RR}$ contains a compact maximal torus $T_{\RR}$. Take a maximal compact subgroup of $M_{\RR}$ containing $T_{\RR}$, and let $\mf{m}_{\RR}= \mf{k}_{\RR} \oplus \mf{p}_{\RR}$ be the associated Cartan decomposition. Then $h \colon \mbb{S} \to M_{\RR}$ will yield an $\Ad(h(\imath))$-polarizable Hodge structure on $\mf{m}$ if $\Ad h(\imath)$ is a Cartan involution of $M_{\RR}$ (or $\mf{m}_{\RR}$). In particular, if we define $h$ as the weight-zero extension of a co-character $\mf{l} \colon \mbb{S}^1 \to T_{\RR}$ satisfying
\begin{align*}
&\text{$\langle \mf{l}, \alpha \rangle \equiv 0 \pmod 4$ for all compact roots $\alpha$;} \\
&\text{$\langle \mf{l}, \beta \rangle \equiv 2 \pmod 4$ for all non-compact roots $\beta$,}
\end{align*}
then $\Ad h(\imath)$, which acts as $\imath^{\langle \mf{l}, \gamma \rangle}$ on the $\gamma$ root space, is trivial on the compact root spaces and is $-1$ on the non-compact root spaces. Thus $\Ad h(\imath)$ is a Cartan involution.\footnote{Explicitly (following IV.B.3 of \cite{griffiths-green-kerr:MTgroups}), the Killing form $B$ is negative-definite on $\mf{k}_0$ and positive-definite on $\mf{p}_0$, so $-B$ gives a polarization in this case.}

Next we must check that we can arrange that $M$ be the full M-T group. Choose any non-trivial $h$ as in the previous part of the proof (i.e. associated to a non-zero element of the polarizable congruence classes), and let $M_1 \subset M$ be the smallest algebraic subgroup (over $\Q$) of $M$ through which all $M(\RR)$-conjugates of $h$ factor. Clearly $M_1$ is a non-trivial normal subgroup of $M$, hence ($M$ is assumed $\Q$-simple) $M_1=M$. But $M_1$ is in fact realized as the M-T group of some conjugate of $h$: indeed, it is the \textit{generic} M-T group for this family, i.e. for all points $h'$ of $M(\RR) \cdot h= M(\RR)/\mr{Stab}(h)$ outside a countable union of closed analytic subspaces, the M-T group of $h'$ is equal to $M_1=M$.\footnote{For this standard fact, see for instance \cite[6.4]{moonen:introMTgroups}; note that the argument only requires a holomorphically varying family of Hodge structures, \textit{not} a variation satisfying Griffiths transversality.}
\endproof
The following extension of the proposition is stated without proof as IV.A.3 of \cite{griffiths-green-kerr:MTgroups}:\footnote{Again, changing `$M$ has an anisotropic maximal torus' to `$M_{\RR}$ has a compact maximal torus.'}
\begin{cor}\label{semisimple}
Let $M/\Q$ be a semi-simple adjoint group. Then $M$ is a M-T group if and only if $M_{\RR}$ contains a compact maximal torus.
\end{cor}
\proof
Note that the proof of Proposition \ref{simpleadjoint} yields infinitely many non-conjugate (over $\Q$) surjections $\mc{H} \onto M'$ for any $\Q$-simple adjoint group $M'$ whose real points contain a compact maximal torus. The result will follow from a standard application of Goursat's lemma for Lie algebras. Let us recall what this says:
\begin{itemize}
\item Suppose that $\mf{h} \subset \mf{m}_1 \oplus \mf{m}_2$ is an inclusion of Lie algebras (over any field), such that $\mf{h}$ surjects onto each $\mf{m}_i$ via the projections $\mf{m}_1 \oplus \mf{m}_2 \xrightarrow{\pi_i} \mf{m}_i$. Then $\mf{h} \cap \ker(\pi_2)$ is an ideal of $\mf{m}_1$ (likewise, switching the roles of $1$ and $2$), and the `graph' of $\mf{h}$ induces an isomorphism of Lie algebras
\[
\mf{m}_1/\left(\mf{h} \cap \ker(\pi_2)\right) \xrightarrow{\sim} \mf{m}_2/\left( \mf{h} \cap \ker(\pi_1)\right).
\]
\end{itemize}
The proof of this fact is immediate. Now, we want to find a surjection $\mc{H} \onto M$; as long as we have a map, we can check surjectivity at the level of Lie algebras. Decomposing the Lie algebra $\mf{m}$ of $M$ as a Lie-algebra direct sum
\[
\mf{m}= \bigoplus_{i=1}^r \mf{m}_i^{\oplus {n_i}},
\]
where the various $\mf{m}_i$ are $\Q$-simple and non-isomorphic, we quickly reduce by Goursat to the case $r=1$, say $M= M_1^{n_1}$ (recall that $M$ is of adjoint type). Now, the group of outer automorphisms (a $\Q$-group scheme) of $\mf{m}_1$ is finite, so in particular, $\Aut(\mf{m}_1)(\Q)/\mr{Inn}(\mf{m}_1)(\Q)$ is finite; choose representatives $\mc{O} \subset \Aut(\mf{m}_1)(\Q)$ of this quotient. We can then (since a morphism $h \colon \mc{H} \to M_1$ is determined by $\mr{Lie}(h)$) find surjections $h_1, \ldots, h_{n_1} \colon \mc{H} \onto M_1$ such that no element of $\mc{O} \cdot \mr{Lie}(h_i)$ is $\mr{Inn}(\mf{m}_1)(\Q)= M_1(\Q)$-conjugate to any element of $\mc{O} \cdot \mr{Lie}(h_j)$ for $i \neq j$. We claim that the map $\oplus \mr{Lie}(h_i) \colon \mf{h} \to \mf{m}_1^{\oplus n_1}$ is surjective. Denote by $\bar{\mf{h}}$ the image, and let $\pi_i \colon \mf{m}_1^{\oplus n_1} \onto \mf{m}_1$ be the projection onto the $i^{th}$ direct factor. By the Lemma to \cite[Theorem 4.4.10]{ribet:AVreal},\footnote{Briefly, we can induct, assuming $\mf{h}$ surjects onto $\oplus_{i=1}^{n_1-1} \mf{m}_1$. Then $\bar{\mf{h}} \cap \ker(\pi_{n_1})$ is an ideal of $\oplus_{i=1}^{n_1-1} \mf{m}_1$, which must be the whole thing since it surjects onto each $\im(\pi_i)$, $1 \leq i \leq n_1-1$.} it suffices to check this surjectivity for each of the pair-wise projections $\mf{h} \to \mf{m}_1^{\oplus n_1} \xrightarrow{(\pi_i, \pi_j)} \mf{m}_1 \oplus \mf{m}_1$. But if inside $\mf{m}_1 \oplus \mf{m}_1$ we had $(\pi_i, \pi_j)(\bar{\mf{h}}) \cap \ker(\pi_i)= (0)$, then $(\mr{Lie}(h_i)(x), \mr{Lie}(h_j)(x))$ would be the graph of an isomorphism $\mf{m}_1 \xrightarrow{\sim} \mf{m}_1$, so that for some automorphism $\alpha$ of $\mf{m}_1$, $\alpha \circ \mr{Lie}(h_i)= \mr{Lie}(h_j)$. By construction of the $h_k$, we have avoided this possibility, so the proof is complete.
\endproof
\begin{rmk}\label{semisimplebis}
The above argument in the semi-simple case is rather more refined than it needs to be. Suppose $M$ is $\Q$-simple adjoint, with a maximal compact torus $T_{\RR}$ inside $M_{\RR}$; the argument of Proposition \ref{simpleadjoint} shows that for \textit{any} non-trivial co-character $\mf{l} \colon \mbb{S}^1 \to T_{\RR}$ satisfying the congruences appearing in that proof, we can find after conjugating (so replacing $T_{\RR}$ with some conjugate torus in $M_{\RR}$) an $h \colon \mbb{S} \to M_{\RR}$ for which $M$ is the M-T group of the induced (polarizable) Hodge structure on $\mf{m}_{\RR}$. In particular, we can do this for infinitely many distinct co-characters $\mf{l}$ satisfying the necessary congruences, say giving rise to Hodge structures $h_i \colon \mbb{S} \to M_{\RR}$.\footnote{These $h_i$ may land in different compact maximal tori of $M_{\RR}$, but when we conjugate them into a common torus $T_{\RR}$ they correspond to different elements of $X_{\bullet}(T)$.} The result is different Hodge structures $\Ad(h_i)$ on $\mf{m}$, from which it follows easily by the above Goursat argument that $\mc{H} \xrightarrow{h_1 \times \cdots \times h_k} M^k$ is surjective.
\end{rmk}
Finally, we complete the proof of Theorem \ref{main}:
\begin{thm}
Let $\tM$ be a connected reductive group over $\Q$, with our running assumptions as in Notation \ref{hypotheses2}. Suppose moreover that there is a surjection $\mc{S} \onto Z^0_{\tM}$, and that, for some Cartan involution of $M^{ad}_{\RR}$, some member $\mu \in X_{\bullet}(T^{ad})$ of the polarizable congruence classes lifts to $\tM_{\CC}$. Then $\tM$ is a quotient of $\mc{H}$, i.e. it is a M-T group.
\end{thm}
\proof
Suppose $\mu \in X_{\bullet}(T^{ad})$ lifts to $\tM_{\CC}$. Then for any even integer $N$ annihilating $X^{\bullet}(\tZ_M)_{tor}$, each element of $\mu + N X_{\bullet}(T^{ad})$ also lifts to $\tM_{\CC}$, and still lies in the polarizable congruence classes. As in Remark \ref{semisimplebis}, any suitably generic element $\mu'$ of $\mu+NX_{\bullet}(T^{ad})$ will yield a Hodge structure with M-T group $M^{ad}$.\footnote{That is, take $\mf{l}= 2\mu \colon \mbb{S}^1 \to T^{ad}_{\RR}$ and, working separately with each $\Q$-simple factor, take suitable $M^{ad}_{\RR}$-conjugates as in the proof of Corollary \ref{semisimple}.} Thus (possibly changing $T^{ad}_{\RR}$), we may assume that we have a Hodge structure $h \colon \mbb{S} \to T^{ad}_{\RR}\subset M^{ad}_{\RR}$ whose M-T group is $M^{ad}$, and whose associated Hodge co-character $\mu_h$ lifts to $\tilde{\mu} \colon \mbb{G}_{m, \CC} \to \tT_{\CC}$. Then $\tilde{h}(z)= \tilde{\mu}(z) \overline{\tilde{\mu}(z)}$ defines a homomorphism $\mbb{S} \to \tM_{\RR}$. Note that $\tilde{h}(\mbb{G}_{m, \RR})$ is contained in $Z^0_{\tM, \RR}$, so that $w_{\tilde{h}}= \tilde{h}|_{\mbb{G}_{m, \RR}}$ is necessarily defined over $\Q$, and if $Z^0_{\tM, \RR}$ contains a copy of $\mbb{G}_{m, \RR}$, we may twist $\tilde{h}$ (i.e. ensure we are in non-zero weight) so that the M-T group $\mr{M}_{\tilde{h}}$ of $\tilde{h}$ contains $\mbb{G}_{m, \Q}$.  $\tilde{h}(\mbb{G}_{m, \RR})$ being central also implies $\Ad \tilde{h}(\imath)$ is an involution; since $Z_{M', \RR}$ is compact and $\Ad h(\imath)$ is a Cartan involution of $M^{ad}_{\RR}$, we conclude that $\Ad \tilde{h}(\imath)$ is a Cartan involution of $M'_{\RR}$. It follows by Lemma \ref{deligne} that any homogeneous (with respect to $w_{\tilde{h}}$) rational representation of $\tM$ is $\Ad \tilde{h}(\imath)$-polarizable, so $\mr{M}_{\tilde{h}}$ is the M-T group of a (polarizable) Hodge structure. Of course, $\mr{M}_{\tilde{h}}$ also contains $M= \tM^{der}$, so we are just left to check that by further twisting $\tilde{h}$ we can arrange $\mr{M}_{\tilde{h}} \supset Z^0_{\tM}$. Let us abusively write $\tilde{h}$ also for the morphism $\tilde{h} \colon \mc{H} \to \tM$; certainly $\tilde{h}(\mc{C}) \subset Z^0_{\tM}$. Let $Z'$ be a $\Q$-torus that is an isogeny complement to $\tilde{h}(\mc{C})$ inside $Z^0_{\tM}$ ($Z'$ exists by Maschke's theorem), so for some integer $N$ we have inclusions of $\Z[\Gal(\Qb/\Q)]$-modules:
\[
X^{\bullet}(Z') \oplus X^{\bullet}(\tilde{h}(\mc{C})) \supset X^{\bullet}(Z^0_{\tM}) \supset N X^{\bullet}(Z') \cong X^{\bullet}(Z').
\] 
Now, by assumption there is a surjection $\mc{S} \onto Z^0_{\tM}$, i.e. an inclusion $X^{\bullet}(Z^0_{\tM}) \subset X^{\bullet}(\mc{S})$. Composing with the above inclusion $X^{\bullet}(Z') \into X^{\bullet}(Z^0_{\tM})$, we obtain a surjection $\gamma \colon \mc{S} \onto Z'$. Twisting $\tilde{h}$ by $\gamma$, we obtain the desired surjection $\tilde{h} \otimes \gamma \colon \mc{H} \onto \tM$.
\endproof

\bibliographystyle{amsalpha}
\bibliography{biblio.bib}

\end{document}